\ifx\LocalElevenPtArticleFormatLoaded\undefined
  \documentclass[a4paper,11pt]{article}     
  \usepackage{array}                        
  \usepackage{theorem}                      
  \usepackage{amsmath,amscd,amssymb}        
  \usepackage{latexsym}                     
  \usepackage[german,english]{babel}        
  \usepackage[applemac]{inputenc}           
                \usepackage{graphicx}
  \usepackage{epsfig}
  \usepackage{xypic}

\fi

\newtheorem{theorem}{Theorem}[section]
\usepackage{theorem}

\newtheorem{lemma}[theorem]{Lemma}
\newtheorem{proposition}[theorem]{Proposition}

\theorembodyfont{\rm}
\newtheorem{remark}[theorem]{Remark}

\newcommand{\CC}{{\mathbb{C}}}

\newcommand{\HH}{{\mathbb{H}}}

\newcommand{\PP}{{\mathbb{P}}}
\newcommand{\QQ}{{\mathbb{Q}}}
\newcommand{\RR}{{\mathbb{R}}}
\newcommand{\ZZ}{{\mathbb{Z}}}

\newcommand{\on}[1]{\operatorname{#1}}
\renewcommand{\phi}{\varphi}
\newenvironment{Proof}{\begin{ProofwCaption}{Proof}}{\end{ProofwCaption}}
\newenvironment{Proof*}[1]{\begin{ProofwCaption}{{#1}}}{\end{ProofwCaption}}
\newenvironment{ProofwCaption}[1]%
  {\addvspace\theorempreskipamount \noindent{\it #1.}\rm}%
  {\qed \par \addvspace\theorempostskipamount}
\newcommand{\qedsymbol}{\mbox{$\Box$}}
\newcommand{\qed}{\quad\qedsymbol}
\begin{document}
\title{The modular form of the Barth-Nieto quintic}

\author{V.~Gritsenko and K.~Hulek}
\date{}
\maketitle

\section {Introduction}
Barth and Nieto showed in \cite{BN} that the quintic threefold
$$
N=\left\{ \sum\limits^5_{i=0} u_i=\sum\limits^5_{i=0}\frac 1{u_i}=0\right\}
\subset
\PP^5.
$$
parametrizes birationally the space of Kummer surfaces associated to abelian
surfaces with a $(1,3)$-polarization and a level $2$ structure. The quintic
$N$ has a smooth model which is a Calabi-Yau threefold. From this Barth and
Nieto deduced that the Siegel modular variety ${\cal A}_3(2)$ parametrizing
abelian surfaces with a $(1,3)$-polarization and a level $2$ structure also
has a smooth model which is a Calabi-Yau manifold
whose Euler number is $80$.
As a consequence this implies that there
is exactly one weight $3$ cusp form (up to a scalar)
with respect to the modular group $\Gamma_3(2)$ which defines the space
${\cal A}_3(2)$. It is a natural question to ask to determine this cusp
form. In
this paper we give the answer to this problem by showing that the modular form
in question is
$\Delta_1^3$
where $\Delta_1$ is a remarkable cusp form of weight $1$
with respect to $\Gamma_3$ (the paramodular group of
a $(1,3)$-polarization) with a character of order $6$.
This function was discovered in \cite{GN2} and it determines
a generalized Lorentzian Kac--Moody superalgebra of Borcherds type.
The form
$\Delta_1$ is  also the automorphic discriminant of the  moduli
space of $K3$ surfaces  whose lattice of transcendental cycles
is  contained in
$U(12)\oplus U(12)\oplus <2>$ (see \cite[Theorem 5.2.1]{GN2}).

One of the main features of $\Delta_{1}$
is that it vanishes precisely along the
diagonal
$${\cal H}_1=\left\{ \tau_2=0\right\}=
\left\{\tau=\left(
\begin{array}{cc}
\tau_1 & 0\\ 0 & \tau_3\end{array}\right),\  \tau_1, \tau_3 \in
\HH_1\right\}\subset \HH_2
$$
which parametrizes split abelian surfaces. The vanishing order along the
diagonal is $1$. In fact it turns out (see Proposition \ref{prop11}) that
there is only a list of four modular forms with respect to the paramodular
group $\Gamma_t$ with this property. All of these modular forms are
particularly interesting (see remarks at the end of Section 1).
This paper is mostly concerned with the geometric
consequences which can be derived from the form $\Delta_1$. In section $2$ we
give a straightforward construction of a smooth, projective
Calabi-Yau model of the variety ${\cal A}_3(2)$, resp. its Voronoi
compactification ${\cal A}^{\ast}_3(2)$
(cf. Theorem \ref{theo27}). Our method is entirely within the
framework of Siegel modular varieties, i.e. it uses the
toroidal compactification (we describe its properties in \ref{theo24}),
but is independent of the arguments of
Barth and Nieto who go via the embeddings of the Kummer surfaces. In
section $3$
we study the modular variety ${\cal A}_3(3)$ where the level $2$ structure is
replaced by a level $3$ structure. Here we can use the existence of the
modular form $\Delta_1$ to show that this space is of general typ and to
construct a minimal model (see Theorems \ref{theo31} and \ref{theo32}).

The first author wants to thank A. Tyurin, who has drawn
his attention to the paper of Barth and Nieto. Both authors are grateful to the
DFG for partial support under grant 436 RUS 17/104/97. The first author also
wants to thank RIMS for its hospitality and support.

\section{Modular forms vanishing along the diagonal}
For an integer $t\ge 1$ the {\em paramodular group} $\Gamma_t$ is the subgroup
$$
\Gamma_t= \left\{
g\in \on {Sp} (4,\QQ); g\in \left(
\begin{array}{cccc}
\ZZ & \ZZ & \ZZ & t\ZZ\\
t\ZZ & \ZZ & t\ZZ & t\ZZ\\
\ZZ & \ZZ & \ZZ & t\ZZ\\
\ZZ & t^{-1}\ZZ & \ZZ & \ZZ
\end{array}
\right)
\right\}
$$
of $\on{Sp}(4,\QQ)$. Its geometric meaning is that the quotient
$$
{\cal A}_t=\Gamma_t\backslash\HH_2
$$
of the Siegel space $\HH_2$ by $\Gamma_t$ is the moduli space of
$(1,t)$-polarized abelian surfaces. Let $\pi:\HH_2\rightarrow {\cal A}_t$
be the quotient map.\\

The diagonal
$${\cal H}_1=\left\{ \tau_2=0\right\}=
\left\{\tau=\left(
\begin{array}{cc}
\tau_1 & 0\\ 0 & \tau_3\end{array}\right),\  \tau_1, \tau_3 \in
\HH_1\right\}\subset \HH_2
$$
resp. its image $\pi({\cal H}_1)$ in ${\cal A}_t$ parametrizes split polarized
abelian surfaces. The surface $\pi({\cal H}_1)$ is a component of the
Humbert surface
$H_1$ of discrimant $1$.
For the theory of Humbert surfaces we refer the reader to
\cite{vdG}, \cite{GH1}. The equality $\pi({\cal H}_1)=H_1$
holds if
and only if
the equation $b^2 \equiv 1 \operatorname{mod} 4t$ has
$\on{mod} 2t$  only
the solution
$\pm 1$.
There is a remarkable series of modular forms (with a
character) whose zero locus in ${\cal A}_t$ consists
exactly of $\pi({\cal H}_1)$.
We list these modular forms in the following table\\

\begin{tabular}{c|c|c|c|c}
form & weight & group & order of character & cusp form\\
\hline
$\Delta_5$ & 5 & $\Gamma_1$ & 2 & yes\\
$\Delta_2$ & 2 & $\Gamma_2$ & 4 & yes\\
$\Delta_1$ & 1 & $\Gamma_3$ & 6 & yes\\
$\Delta_{1/2}$ & 1/2 & $\Gamma_4$ & 8 &no
\end{tabular}\\

\vspace{5mm}

\noindent
These forms were discussed in some detail in \cite{GN2}. They arise from Jacobi
forms by arithmetical lifting and are denominator
functions of generalized Kac-Moody superalgebras. These functions  vanish on
${\cal H}_1$ (and its $\Gamma_t$-translates) of order 1 and nowhere else. Such
forms are very special, in fact we have

\begin{proposition}\label{prop11}
 If $F$ is a modular form  of integral
(or half-integral) weight with a character (or a multiplier system)
with respect to $\Gamma_t$ such that $F$ vanishes exactly on the
$\Gamma_t$-translates of
${\cal H}_1$ and with vanishing order $1$, then
$t=1$, $2$, $3$ or  $4$ and $F$ is equal, up to a constant,
to $\Delta_5$, $\Delta_2$, $\Delta_1$ or $\Delta_{1/2}$.
\end{proposition}

\begin{Proof}
In the theory of automorphic forms it is sometimes more
natural to deal with a group conjugated to $\Gamma_t$, namely
$$
\Gamma'_t=I_t\Gamma_t I_t^{-1},
$$
where $I_t=\hbox{diag}(1,t^{-1},1,t)$. This is again a subgroup of
$\on{Sp}(4,\QQ)$.
Note that ${\cal A}'_t=\Gamma'_t\backslash\HH_2\cong
\Gamma_t\backslash \HH_2={\cal A}_t$. In order to avoid confusion we shall
denote all objects which refer to the group $\Gamma'_t$ by a $'$.
Let $F$ be a modular form of weight $k$ with respect to $\Gamma'_t$
and $d\ne t$ be an integer.
We shall use the following operator of multiplicative symmetrisation
\begin{eqnarray}
[F]_{d}=
\prod\limits_{ M\in \Gamma'_{t}\cap \Gamma'_{d}\setminus \Gamma'_{d}}
F|_k M
\end{eqnarray}
where
$(F|_kM)(Z):=\hbox{det\,}(C\tau+D)^{-k}F(M<Z>)$
($M=\left(\smallmatrix A&B\\C&D \endsmallmatrix\right)$)
is the standard slash-operator.
(Compare (1) with the operators of symmetrisation studied
in \cite{G2} and \cite[\S 3]{GN2}).
It is clear that $[F]_{d}$ is a modular form with respect to $\Gamma'_{d}$.

A rational quadratic divisor ${\cal H}'_l$ with respect to $\Gamma'_{t}$
is defined by
$$
{\cal H}'_{l}=
\{\pmatrix \tau_1&\tau_2\\ \tau_2&\tau_3 \endpmatrix\in \HH_2\ |\
f({\tau_2}^2-\tau_1\tau_3)+c\tau_3+b\tau_2+a\tau_1+e=0\,\}.
$$
where $l=(e,a,b,c,f)\in \ZZ^3\times (t\ZZ)^2$ such that
$(e,a,b,\frac c{t}, \frac{f}{t})=1$.
The integer
$D(l)=b^2-4ac-4ef$ is called the discriminant of ${\cal H}'_{l}$.
The surface $H'_l$ in ${\cal A}'_t$ is defined as the image of ${\cal H}'_l$
under the natural projection map $\pi'_t: \HH_2\to {\cal A}'_t$. It is a
component of the Humbert surface $H'_{D(l)}$ and by abuse of notation we shall
also sometimes refer to $H'_l$ as an Humbert surface.

Let $d$ be a divisor of $t$. Then ${\cal H}'_l$ is a rational quadratic
divisor of the same discriminant $D(l)$ with respect to $\Gamma'_{d}$,
if $(e,a,b,\frac c{d}, \frac{f}{d})=1$. For example, this is the case
if $(D(l),t)=1$.
Let us assume, that
$\hbox{div}_{{\cal A}'_t} F=m H'_l$
and $(D(l),t)=1$.
It follows from the consideration above that
$
\hbox{div}_{{\cal A}'_{d}} ([F]_d)
$
is a sum of some irreducible components of the
Humbert surface of discriminant $D(l)$.

We are interested in the case $d=1$ and $H'_l=H'_1=H_1$.
Note that the Humbert surface $H_1$ is irreducible in ${\cal A}_1$. We set
${}_0\Gamma'_t=\Gamma_1\cap \Gamma'_t$.
By standard arguments we see that
for an element $g\in \Gamma_1$ the class ${}_0\Gamma'_t g$
is determined by the last line of $g$  considered
as an element of $\PP^3(\ZZ_t^*)$.
More precisely,
$$
|{}_0\Gamma'_t \setminus \Gamma_1| =
|\ZZ_t^*\setminus \{(a,b,c,d) \hbox{ mod }t,\  (a,b,c,d)=1\}|.
$$
Therefore
$$
[\Gamma_1:{}_0\Gamma'_t]=
\phi(t)^{-1}t^4\prod_{p|t}(1-p^{-4})=t^3\prod_{p|t}(1+p^{-1})(1+p^{-2}).
$$
Here $\phi$ as usual denotes the Euler $\phi$-function.
The weight of $[F]_1$ equals the weight of $F$ multiplied by the
above index.

The order   of zero of $[F]_1$ along the Humbert surface $H_1$
is equal to the number of left cosets ${}_0\Gamma'_t M$ in
$\Gamma_1$ such that there exists an element
$\gamma \in \Gamma'_t$ with $\gamma^{-1}M<{\cal H}_1>={\cal H}_1$.
By \cite{Fk} the stabilizer of ${\cal H}_1$ in $\operatorname{Sp}_4(\RR)$
is the group generated by
$$
\operatorname{SL}_2(\mathbb R)\times \operatorname{SL}_2(\mathbb R) \cong
\biggr\{\left(
\smallmatrix
a&0  &b&0\\
0&a_1&0&b_1\\
c&0  &d&0\\
0&c_1&0&d_1
\endsmallmatrix\right)\in \operatorname{Sp}(\mathbb R)\biggl\}
$$
and the involution
\begin{eqnarray}
\left(\smallmatrix
0&1&0&0\\
1&0&0&0\\
0&0&0&1\\
0&0&1&0
\endsmallmatrix\right)
 <\pmatrix \tau_1&\tau_2\\ \tau_2&\tau_3\endpmatrix>=
\pmatrix \tau_3&\tau_2\\ \tau_2&\tau_1\endpmatrix.
\end{eqnarray}
It follows that  $M<{\cal H}_1>$ coincides with
a $\Gamma_t'$-translate of ${\cal H}_1$ if and only if
${}_0\Gamma'_t M$ contains an element with the last line
$(0,*,0,*)$ or $(*,0,*,0)$ $\hbox{mod }t$.
The number of such classes is equal to
$$
2 \phi(t)^{-1}t^2 \prod_{p|t}(1-p^{-2})=
2t \prod_{p|t}(1+p^{-1}).
$$

As a result we have proved the following: if $F$ is of weight $k$ and its
divisor in $\HH_2$ is exactly the $\Gamma'_t$-orbit of $m{\cal H}_1$,
then
$[F]_1$ has weight
$k\, t^3\prod_{p|t}(1+p^{-1})(1+p^{-2})$ and divisor
$$
\hbox{div}_{{\cal A}_1} \bigl([F]_1\bigr)=
2mt \prod_{p|t}(1+p^{-1}) \,H_1.
$$
It is well known (see \cite{F1}), that for the Siegel modular group
$\Gamma_1$ there is a cusp form $\Delta_5$ of weight $5$
with a character of order $2$ whose divisor is $H_1$.
Using the Koecher
principle we conclude that
$[F]_1$ is a power of $\Delta_5$. Thus
\begin{eqnarray}
k t^3\prod_{p|t}(1+p^{-1})(1+p^{-2})=10m t \prod_{p|t}(1+p^{-1}).
\end{eqnarray}
If $m=1$, then we have
$k t^2\prod_{p|t}(1+p^{-2})=10$.
Since $t^2$ must be smaller than $20$, there are only four numerical
possibilities, namely
$(t,k)=(1,5)$, $(2,2)$, $(3,1)$, $(4, \frac 1{2})$.
Again by Koecher's principle the corresponding
modular forms are, if they exist, unique up to a
scalar, hence it remains to recall the
existence of the modular forms in question.
We have already used the famous
modular form $\Delta_5$, which is the product
of all even Siegel theta-constants (see \cite{F1}).
For $t=4$ we consider the even Siegel theta-constant with characteristic
$\frac 1{2}
\left(\smallmatrix 1&1\\1&1\endsmallmatrix\right)$:
\begin{eqnarray}
\Delta_{1/2}(Z)&=
\frac{1}2\sum\limits_{ n,\,m\in \ZZ}
\,\bigl(\frac{-4}{n}\bigr) \bigl(\frac{-4}{m}\bigr)
\exp{\bigl(\pi i(\frac{n^2}4\tau_1+nm\tau_2+m^2\tau_3)\bigl)}
\nonumber
\\
{}&=
\sum\limits_{m>0} \bigl(\frac{-4}{m}\bigr)
\,\vartheta(\tau_1, m\tau_2)\exp{(\pi i m^2 \tau_3)}.
\nonumber
\end{eqnarray}
Here $\vartheta(\tau_1,\tau_2)$ is the Jacobi
theta-constant with characteristic $(\frac 1{2},\frac 1{2})$,
which is a Jacobi modular form of index $\frac 1{2}$.
This representation implies (see \cite[Theorem 1.11]{GN2}) that
$\Delta_{1/2}$ is a modular form with respect to $\Gamma'_4$
since $\Gamma'_4$ is generated by the Jacobi group and the
involution
$$
\pmatrix \tau_1&\tau_2\\ \tau_2&\tau_3\endpmatrix \mapsto
\pmatrix {\tau_3}/4 &\tau_2\\ \tau_2& 4\tau_1\endpmatrix.
$$
The cusp forms $\Delta_2$ and $\Delta_1$ for $\Gamma'_2$
and $\Gamma'_3$ can be defined as the arithmetic lifting of the Jacobi
forms
$\phi_{1}(\tau_1,\tau_2)=\eta(\tau_1)^3\vartheta(\tau_1,\tau_2)$
and
$\phi_{1}(\tau_1,\tau_2)=\eta(\tau_1)\vartheta(\tau_1,\tau_2)$
of index $\frac{1}2$ respectively (see \cite[Theorem 1.12]{GN2}):
\begin{eqnarray}
\Delta_2(Z)&=
\sum\limits_{\scriptsize
\begin{array}{c} m\equiv 1\,mod\,4\\
 m>0\end{array}}
m\sum\limits _{\scriptsize\begin{array}{c}ad=m\\  b\,mod\,d\end{array}}
d^{-2}\bigl(\frac {-4}a\bigr)\phi_1\left(\frac{a\tau_1
+4b}d,\,a\tau_2\right)
\exp{(\pi i m\tau_3)}
\nonumber
\end{eqnarray}
\begin{eqnarray}
\Delta_1(Z)&=
\sum\limits_{\scriptsize\begin{array}{c} m\equiv 1\,mod\,6\\
 m>0\end{array}}
\sum\limits_{\scriptsize\begin{array}{c} ad=m\\  b\,mod\,d\end{array}}
d^{-1}\bigl(\frac {-4}a\bigr)\phi_2(\frac{a\tau_1 +6b}d,\,a\tau_2)
\exp{(\pi i m\tau_3)}.
\nonumber
\end{eqnarray}
$\Delta_{1/2}$, $\Delta_1$, $\Delta_2$ vanish along ${\cal H}_1$,
because $\vartheta(\tau_1, 0)\equiv 0$.
The consideration with $[F]_1$ above implies that the order of zero
along the diagonal for each modular form is one and that $H_1$
is its full divisor.
(In \cite{GN2} this fact was proved using the Borcherds lifting.)
Representations of the modular forms
$\Delta_1$, $\Delta_2$ as lifting give us  elementary formulae
for the Fourier coefficients of these functions
(see \cite[Example 1.14]{GN2})
which imply, for example, that $\Delta_1$ and $\Delta_2$ are cusp forms:
$$
\hskip -2pt\Delta_1(Z)=\hskip -4pt\sum_{M\ge 1}\hskip -7pt
\sum\limits_{\scriptsize
\begin{array}{c}
 n,\,m >0,\,l\in \ZZ\\
n,\,m\equiv 1\,mod\,6\\
4nm-3l^2=M^2
\end{array}}\hskip -8pt
\biggl(\frac{-4}{l}\biggr)\hskip -2pt
\biggl(\frac{12}{M}\biggr)\hskip -6pt
\sum\limits_{a|(n,l,m)}
\hskip -3pt\biggl(\frac{6}{a}\biggr)
\exp{\bigl(\pi i(\frac{n}3\tau_1+l\tau_2+m\tau_3)\bigl)}
$$
and
$$
\hskip -2pt\Delta_2(Z)=\hskip -4pt\sum_{N\ge 1}\hskip -6pt
\sum\limits_{\scriptsize
\begin{array}{c}
 n,\,m >0,\,l\in \ZZ\\
n,\,m\equiv 1\,mod\,4\\
2nm-l^2=N^2
\end{array}}\hskip -8pt
N\biggl(\frac {-4}{Nl}\biggr)\hskip -4pt
\sum_{a\,|\,(n,l,m)} \hskip -3pt\biggl(\frac {-4}{a}\biggr)
\exp{\bigl(\pi i(\frac{n}2\tau_1+l\tau_2+m\tau_3)\bigl)}.
$$

\hfill\end{Proof}

The symmetrisation (1) gives us some formulae for
$\Delta_5$ and $\Delta_2$ in terms of the theta-constant $\Delta_{1/2}$.

{\bf 1. $\Delta_5$ in terms of $\Delta_{1/2}$}. First  we have seen
in the above proof that
$$
\Delta_5^2= (const)\prod\limits_{ M \in \Gamma'_4\cap\Gamma_1\setminus
\Gamma_1}
\Delta_{1/2}|_{\frac 1{2}} M.
$$
One can get a simpler representation (see \cite[Theorem 1.11]{GN2})
considering  a subgroup $\Gamma_{1,2}$
of $\Gamma_1$ conjugated to $\Gamma'_4$, namely
$$
\Gamma_{1,2}=\hbox{diag}(1,2,1,2^{-1})\,\Gamma'_4\,
\hbox{diag}(1,2^{-1},1,2).
$$
The modular form
$\tilde{\Delta}_{1/2}(Z)={\Delta}_{1/2}
(\left(\smallmatrix \tau_1& \tau_2/2\\ \tau_2/2&\tau_3/4
\endsmallmatrix\right))$
is $\Gamma_{1,2}^+$-modular form, where $\Gamma_{1,2}^+$
is the double extension of $\Gamma_{1,2}$ defined by the
involution (2).
It is easy to check that
$[\Gamma_1:\Gamma_{1,2}^+]=10$.
Thus
$$
\Delta_5(Z)= (const)\prod\limits_
{M \in \Gamma_{1,2}^+ \setminus \Gamma_1}
\tilde{\Delta}_{1/2}|_{\frac 1{2}} M (Z).
$$
This is a new variant of the classical representation of $\Delta_5(Z)$
as  product of $10$ theta-constants.

{\bf 2. $\Delta_2$ in terms of $\Delta_{1/2}$}.
Here we can use the symmetrisation for the pair  $(\Gamma'_2,
\Gamma'_4)$.
We fix a system of representatives
${}_0\Gamma'_4\setminus {}_0\Gamma'_2=
\{M_1,\dots , M_8\}$.
It is easy to see that
$$
{}_0\Gamma'_2\setminus \Gamma'_2=
\biggl\{E_4,\
S_2=\left(\smallmatrix
1&0&0&0\\
0&1&0&2^{-1}\\
0&0&1&0\\
0&0&0&1\endsmallmatrix\right), \
J_2=\left(\smallmatrix
0&\hskip 1pt0&-1&0\\
0&0&\hphantom{-}0&-2^{-1}\\
1&0&\hphantom{-}0&0\\
0&2&\hphantom{-}0&0\endsmallmatrix\right)
\biggr\}.
$$
Since $S_2\in \Gamma'_4$ it follows that
$$
\Gamma'_2\cap\Gamma'_4 \setminus \Gamma'_2=
\{M_1,\dots , M_8,\, M_1J_2,\dots , M_8J_2\}.
$$
As in the proof of Proposition 1.1, it follows  that
$$
[\Delta_{1/2}]_2=\prod\limits_
{M \in \Gamma'_2\cap\Gamma'_4 \setminus \Gamma'_2}
\Delta_{1/2}|_{\frac 1{2}} M= (const)\, \Delta_2^4.
$$
We note that $\Delta_2^4$ is the cusp form of minimal weight
with trivial character for $\Gamma_2$ (see \cite{F3}, \cite{G2}).\\

The modular forms considered above have some applications to
algebraic geometry and physics.
It is known  (see \cite[Theorem 5.2.1]{GN2} and  \cite{GN3}) that they
are the automorphic discriminant of a moduli space of
$K3$ surfaces with the lattice of transcendental cycles
of type $U(n)^2\oplus <2>$, where
$U(n)=\pmatrix 0&-n\\-n&0\endpmatrix$ and
$n=1$ for $\Delta_5$,
$n=8$ for $\Delta_2$,
$n=12$ for $\Delta_1$ and
$n=16$ for $\Delta_{1/2}$.
In \cite{GN1}--\cite{GN2} it was proved that the  three cusp forms
considered above
determine the first members of the main series $A_I$--$A_{III}$ of
generalized Lorentzian Kac--Moody superalgebras of rank $3$ and they
have an interesting infinite product expansion.
In the case of $\Delta_1$ we have the following formula
(see \cite[Theorem 2.6]{GN2}):
$$
\Delta_1(Z)=
q^{\frac {1}6}r^{\frac {1}2}s^{\frac {1}2}
\prod
\limits_{\scriptsize
\begin{array}{c}
 n\ge 0,\,m \ge 0,\,l\in \ZZ\\
(l<0\ if\  n=m=0)
\end{array}}
\bigl(1-q^n r^l s^{3m}\bigr)^{f(nm,l)}
$$
where  $q=e^{2\pi i \tau_1}$, $r=e^{2\pi i \tau_2}$,
$s=e^{2\pi i \tau_3}$ and
$$
\sum_{n\ge 0,\,l}f(n,l)\,q^nr^l=
\biggl(\frac {\vartheta(\tau_1,2\tau_2)}{\vartheta(\tau_1,\tau_2)}\biggr)^2
$$
$$
=r^{-1}
\biggl(\prod_{n\ge 1}(1+q^{n-1}r)(1+q^{n}r^{-1})(1-q^{2n-1}r^2)
(1-q^{2n-1}r^{-2})\biggr)^2.
$$
In physics  the modular  form $\Delta_5$ appears in the two-loops vacuum
amplitude  of bosonic strings (see, for example, \cite{BK}, \cite{M1}).
We remark also that $\Delta_2$ and $\Delta_5$ are related to
the perturbative prepotential and the perturbative Wilsonian gravitational
coupling of some four parameter $D=4$, $N=2$ string models
(see \cite{K1}, \cite{CCL}, \cite{C}).
Moreover in \cite{DVV} it was shown that
$\Delta_5^{-2}$  can be interpreted as the second quantized
elliptic genus of a K3 surface (see also \cite{K2}, \cite{M2}).

In this note we discuss the geometric significance of the form $\Delta_1$ which
is very closely connected with the Barth-Nieto quintic \cite{BN} and we shall
briefly comment on the forms $\Delta_2$ and $\Delta_5$.

\section{The moduli space ${\cal A}_3(2)$}

Let $\Gamma_3(2)$ be the subgroup of the paramodular group $\Gamma_3$ which
defines the moduli space ${\cal A}_3(2)=\Gamma_3(2)\backslash \HH_2$ of
$(1,3)$-polarized abelian surfaces with a level 2 structure. Then
$$
\Gamma_3(2)=\left\{
g\in \Gamma_3; g-{\bf 1} \in\left(
\begin{array}{cccc}
2\ZZ & 2\ZZ & 2\ZZ  & 6\ZZ\\
6\ZZ & 2\ZZ & 6\ZZ  & 6\ZZ\\
2\ZZ & 2\ZZ & 2\ZZ  & 6\ZZ\\
2\ZZ & \frac23\ZZ & 2\ZZ  & 2\ZZ\\
\end{array}
\right)
\right\}.
$$
The group $\Gamma_3$ is conjugate via $\on{diag}(1,1,1,3)$ to the
symplectic group
$\on{Sp}(\Lambda_3,\ZZ)$ where $\Lambda_3$ is the following symplectic form
$$
\left(
\begin{array}{rrcc}
0 & 0 & 1 & 0\\
0 & 0 & 0 & 3\\
-1& 0 & 0 & 0\\
0 &-3 & 0 & 0
\end{array}
\right).
$$
Under this
isomorphism $\Gamma_3(2)$ is identified with the group
$\on{Sp}^{(2)}(\Lambda_3,\ZZ)$ consisting of all elemets $g\in
\on{Sp}(\Lambda_3,\ZZ)$ with $g={\bf 1}\mbox{ mod }2$.

\begin{lemma}\label{lem2.1}
The modular form $\Delta^3_1$ is a weight 3 cusp form with respect to
$\Gamma_3(2)$.
\end{lemma}
\begin{Proof}
The form $\Delta_1$ is a cusp form with respect to $\Gamma_3$ with a character
$\chi_6$ of order 6 \cite{GN2}. Hence $\Delta^3_1$ has a character
$\chi_2=\chi_6^3$ of order 2. It follows from \cite[Theorem 2.1]{GH2} that
there
is exactly one such character and that this character arises in the
following way

$$
\diagram
1 \rto & \Gamma_3(2)\rto &\Gamma_3\drto^{\chi_2}\rto &
\on{Sp}(4,\ZZ/2)\cong S_6
\dto^{\mbox{sign}}\rto &  1.\\
        &        &   & \{\pm 1\} &
\enddiagram
$$
In particular $\chi_2|_{\Gamma_3(2)}\equiv 1$ and hence $\Delta_1^3$ is a
modular
form with respect to $\Gamma_3(2)$.
\hfill\end{Proof}

\begin{remark}\label{rem22}
It follows from the results of Barth and Nieto \cite{BN} that ${\cal A}_3(2)$
has a smooth projective model $\tilde{\cal A}_3(2)$ which is a Calabi-Yau
3-fold. By Freitag's extension result the space $S_3(\Gamma_3(2))$ of weight 3
cusp forms is isomorphic to $\Gamma(\tilde{\cal A}_3(2), K_{\tilde{\cal
A}_3(2)})$ for every smooth projective model. It follows that $\Delta^3_1$ is
the unique weight 3 cusp form with respect to $\Gamma_3(2)$. We shall, however,
not use the result of Barth and Nieto in what follows. In fact our methods will
allow us to construct a smooth projective Calabi-Yau model of ${\cal A}_3(2)$
in an easy way.
\end{remark}

For what follows we have to determine the singularities of the moduli space
${\cal A}_3(2)$ and of a suitable toroidal compactification. An important role
is played by the involution
$$
I:=\left(
\begin{array}{cccc}
1&&&\\
&-1&&\\
&&1&\\
&&&-1
\end{array}\right)
$$
Note that $I\in\Gamma_3(2)$ and that
$$
\mbox{Fix } I={\cal H}_1.
$$
\begin{lemma}\label{lem23}
Up to conjugation with elements in $\Gamma_3$ the only elements of finite order
in $\Gamma_3(2)$ are $\pm{\bf 1}, \pm I$.
\end{lemma}

\begin{Proof}
We work in $\on{Sp}(\Lambda_3,\ZZ)$. Then every element $g\in
\on{Sp}^{(2)}(\Lambda_3,\ZZ)$ is of the form
$$
g=\left(
\begin{array}{cc}
{\bf 1} + A  &  B\\
C & {\bf 1} + D
\end{array}
\right)\mbox{ with } A, B, C, D \equiv 0\on{ mod } 2.
$$
Hence
$$
g^2=\left(
\begin{array}{cc}
{\bf 1} + 2A + A^2 +BC & 2B + AB +BD\\
2C + CA + DC           & {\bf 1} + CB + 2 D + D^2
\end{array}
\right)\equiv{\bf 1 }\on{ mod } 4.
$$
This shows that $g^2={\bf 1}$ and hence $g$ is an involution. By  a result of
Brasch \cite[Folgerung 2.9]{Br} the only involutions in $\Gamma_3$ are
$$
\pm{\bf 1}, \pm I, \pm\left(
\begin{array}{rcrr}
-1 & 0 & 0 & 0\\
-3 & 1 & 0 & 0\\
0  & 0 &-1 &-3\\
0  & 0 & 0 & 1
\end{array}
\right).
$$
Since the last involution is not in $\Gamma_3(2)$, the claim follows.
\hfill\end{Proof}

We consider the toroidal compactification ${\cal A}^{\ast}_3(2)$ which belongs
to the second Voronoi decomposition.  We shall refer to this
compactification either as {\em Voronoi} or as in \cite{HKW2} as
{\em Igusa compactification} of ${\cal A}_3(2)$.
\begin{theorem}\label{theo24}
The variety ${\cal A}^{\ast}_3(2)$ is Gorenstein. It has exactly $15$
isolated
singularities of type $V_{\frac13(1,1,1)}$.
\end{theorem}
\begin{Proof}
It follows from Lemma \ref{lem23} that the branch locus of the projection map
$\HH_2\rightarrow{\cal A}_3(2)$ are the translates of ${\cal H}_1$. Since $I$
acts locally like a reflection the space ${\cal A}_3(2)$ is smooth. It
remains to determine the singularities on the boundary. Here we can proceed
along the same lines as in \cite{HKW1} and \cite{Br}. We shall first treat the
corank 1 boundary components $D(l)$. Up to the action of
$\Gamma_3/\Gamma_3(2)\cong S_6$ there are two types, namely $D(l_0)$ and
$D(l_{(0,1)})$ where $l_0=(0,0,1,0)$ and $l_{(0,1)}=(0,0,0,1)$. Here we shall
treat the boundary surface $D(l_0)$ in detail, the other boundary surface
$D(l_{(0,1)})$ can be treated in exactly the same way. The parabolic subgroup
associated to $l_0$ is
$$\begin{array}{l}
P_{l_0}(\Gamma_3(2))=\\[3mm]
=\left\{\left(
\begin{array}{cccc}
\varepsilon & m & q & 3n\\
0 & a & * &3b\\
0 & 0 & \varepsilon & 0\\
0 & c/3 & * & d
\end{array}
\right) ; m, n, q, \in 2\ZZ, \varepsilon=\pm 1,\left(
\begin{array}{cc}
a & b\\
c & d
\end{array}\right)
\in \Gamma_1(2)\right\}.
\end{array}
$$
Here the entries $*$ are determined by the condition that the
matrix is symplectic. This description follows from \cite[Proposition
I.3.87]{HKW2} and the  description of the group
$\Gamma_3(2)$. Dividing out by the rank $1$ lattice $P'_{l_0}(\Gamma_3(2))$
which
is given by the elements with $\varepsilon=1, m=n=0, \left(
\begin{array}{cc}
a & b\\ c & d
\end{array}\right)={\bf 1}$, we
obtain the quotient group
$P''_{l_0}(\Gamma_3(2))=P_{l_0}(\Gamma_3(2))/P'_{l_0}(\Gamma_3(2))$ which can
be identified with the following matrix group:
$$
P''_{l_0}(\Gamma_3(2))\cong\left\{\left(
\begin{array}{ccc}
1 & \varepsilon m & \varepsilon n\\
0 & \varepsilon a & \varepsilon b\\
0 & \varepsilon c & \varepsilon d
\end{array}
\right)\right.;\left(
\begin{array}{cc}
a & b\\ c & d
\end{array}\right)\in\Gamma_1(2),m,n\in 2\ZZ,\varepsilon=\pm 1\}.
$$
This acts on $\CC^{\ast} \times \CC \times \HH_1$ with coordinates
$t_1=e^{2\pi i \tau_1/2}, \tau_2, \tau_3$ as follows:
$$
\begin{array}{rcl}
t_1'& =& t_1 e^{\pi i\varepsilon [m\tau_2-\tau_2'(\frac c3\tau_2+c\varepsilon
n-d\varepsilon m)]}\\
\tau_1'& =& (\varepsilon\tau_2+m\tau_3+3n)(\frac c3\tau_3+d)^{-1}\\
\tau_3'& = &3(\frac a3 \tau_3+b)(\frac c3\tau_3+d)^{-1}.
\end{array}
$$We have to find all points $P=(t_1,\tau_2,\tau_3)$ with $g(P)=g$ for some
$1\neq g\in P''_{l_0}(\Gamma_3(2))$. Invariance of the third component implies
that $\tau_3/3$ is a fixed point for $\left(\begin{array}{cc}
a & b\\ c & d \end{array}\right)$. Hence $\left(\begin{array}{cc}
a & b\\ c & d \end{array}\right)={\bf 1}_2$ or $\left(\begin{array}{cc}
a & b\\ c & d \end{array}\right)=-{\bf 1}_2$. In the
first case we have the possibilities $\varepsilon=1$ or $-1$. If
$\varepsilon=1$
then
$g=1$. If $\varepsilon=-1$ then
$$
g=\left(
\begin{array}{ccc}
1 & m & n\\
0 & -1 & 0\\
0 & 0 & -1
\end{array}\right)
$$
with fixed locus
$$
\tau_2=\frac 12(m\tau_3+3n).
$$
Since the involution $g$ acts locally like a reflection in the neighbourhood of
this point this leads to smooth points in ${\cal A}_3(2)$. (Note that $g$ is
induced by the involution
$$
\left(
\begin{array}{cccc}
1  &  m  &  0  &  3n\\
0  & -1  &  0  &  0\\
0  &  0  &  1  &  0\\
0  &  0  &  0  & -1
\end{array}
\right)
\in\Gamma_3(2).
$$
This involution is conjugate to $I$ and hence the curve given by
$\tau_2=(m\tau_3+n)/2$ is the intersection of a translate of ${\cal H}_1$ with
the boundary component $D(l_0)$). If $\left(\begin{array}{cc}
a & b\\ c & d \end{array}\right)=-{\bf 1}_2$ the case $\varepsilon=-1$ gives
$g=1$ and $\varepsilon=1$ again leads to the curve $\tau_2=(m\tau_3+n)/2$.\\

It remains to consider the corank 2 boundary components $E(h)$.
Since all of
these are equivalent under $\Gamma_3/\Gamma_3(2)\cong S_6$ is suffices to
consider the case $h=(0,0,1,0)\wedge(0,0,0,1)$. Here
$$
P(h)=\left\{\left(
\begin{array}{cc}
A  &  0\\
0  & {^t A^{-1}}
\end{array}
\right);
A\in \on{GL}(2,\ZZ), A-{\bf 1}_2\in\left(
\begin{array}{cc}
2\ZZ & 2\ZZ\\
6\ZZ & 2\ZZ
\end{array}
\right)\right\}\cdot P'(h)
$$
where $P'(h)$ is the lattice
$$
P'(h)=\left\{\left(
\begin{array}{cc}
{\bf 1}_2  &  B\\
0  & {\bf 1}_2
\end{array}
\right);
B={^t B}, B\in\left(
\begin{array}{cc}
2\ZZ & 6\ZZ\\
6\ZZ & 6\ZZ
\end{array}
\right)\right\}.
$$
Let $N\subset \mbox{Sym}(2,\ZZ)$ be the lattice generated by the matrices $B$
and let $\Sigma_N\subset N_{\RR}$ be the fan induced by the Legendre
decomposition (which is here equal to the second Voronoi decomposition). Then
we have the partial quotient
$$
\begin{array}{rcl}
e(h):\quad \HH_2 & \rightarrow & T_n\cong(\CC^{\ast})^3\\
\left(
\begin{array}{cc}
\tau_1 & \tau_2\\
\tau_2 & \tau_3
\end{array}
\right)
 &\mapsto & \left(e^{2\pi i\tau_1/2}, e^{2\pi i\tau_2/6},
e^{2\pi i\tau_3/6}\right).
\end{array}
$$
The fan $\Sigma_N$ defines a torus embedding $T_N\subset X_{\Sigma_N}=:X$.
Here, however, the situation differs from \cite{HKW1}. The main  difference is
that $N$ is not equivalent to the lattice $\mbox{Sym}(2,\ZZ)$. We can,
however, compare the torus embedding $T_N\subset X_{\Sigma_N}$ with the
standard
torus embedding $T_{\on{Sym}(2,\ZZ)}\subset X_{\Sigma}$. To do this let $N'$ be
the lattice spanned by matrices $B'$ with
$$
B'={^tB'};\quad B'\in\left(
\begin{array}{cc}
6\ZZ & 6\ZZ\\
6\ZZ & 6\ZZ
\end{array}\right).
$$
Then $N'$ is equivalent to $\mbox{Sym}(2,\ZZ)$ and $N/N'\cong\ZZ/3$. This
defines a commutative diagram
$$
\diagram
T_{N'}\dto & \subset &  X'=X_{\Sigma_{N'}} \dto \\
T_N & \subset & X=X_{\Sigma_N}
\enddiagram
$$
where the vertical maps are quotients by the cyclic group $N/N'$. The variety
$X'$ is covered by affine sets $X'_{\sigma}\cong\CC^3$ and hence smooth.
This is
no longer the case for $X$. In the analogous situation without a level $2$
structure Brasch \cite[Satz(III.5.21)]{Br} has shown that $X_{\Sigma_N}$
has, up
to the action of $P''(h)$, exactly one singularity $P$ which is of
type $V_{\frac13(1,1,1)}$.
This singularity arises as follows: Let $\sigma_1$ be the cone
$$
\sigma_1=\RR_{\ge 0}\left(
\begin{array}{cc}
1 & 1\\
1 & 1
\end{array}\right)+\RR_{\ge 0}\left(
\begin{array}{cc}
1 & 2\\
2 & 4
\end{array}\right)+R_{\ge 0}\left(
\begin{array}{cc}
0 & 0\\
0 & 1
\end{array}\right).
$$
Then $X'_{\sigma_1}\cong\CC^3$ and $N/N'$ acts on $X'_{\sigma_1}$ by
$\langle \mbox{diag}(\rho)\rangle$ where $\rho=e^{2\pi i/3}$.
Hence the
origin gives rise to a
singularity $P$ of type $V_{\frac 13(1,1,1)}$. A similar analysis
applies in the presence of a level $2$ structure. A straightforward calculation
(see also \cite[Hilfssatz(III.5.23]{Br})
shows that the stabilizer of $P$ in $P''(h)$ is trivial
if we work with a level $2$ structure, resp. isomorphic to the symmetric
group $S_3$ if we work without a level structure. This implies that the
singularities of ${\cal A}^{\ast}_3(2)$ are of type $V_{\frac13(1,1,1)}$.
To count the number of singularities recall that the number of corank $2$
boundary components is $15$ (see e.g. \cite{Fr}). It remains to show that
every boundary component $E(h)$ has exactly one singularity. To count this
number we consider the groups
$$
\Gamma_0(3)=\left\{\left(
\begin{array}{cc}
a & b\\
c & d
\end{array}\right)\in\on{SL}(2,\ZZ); c\equiv 0\mbox{ mod } 3\right\}.
$$
and
$$
\Gamma_0^{(2)}(3)=\left\{\left(
\begin{array}{cc}
a & b\\
c & d
\end{array}\right)\in\Gamma_0(3); \left(
\begin{array}{cc}
a & b\\
c & d
\end{array}\right)\equiv{\bf 1}\mbox{ mod } 2\right\}.
$$
The quotient $\Gamma_0(3)/\Gamma_0^{(2)}(3)$ is isomorphic to the symmetric
group $S_3$ and this is exactly the stabilizer of the cone $\sigma_1$,
resp. the point $P$. This gives the claim.
\hfill\end{Proof}

We have already introduced the Humbert surface $H_1\subset{\cal A}_3$. The
number of components of $H_1$ in ${\cal A}_3$ is given by $\#\{\pm
\  b\mbox{ mod }6; b^2 \equiv 1 \mbox{ mod }12\}=1$ (see \cite{GH2}). Hence
$H_1=\pi({\cal H}_1)$. Recall that $H_1$ parametrizes
split polarized abelian surfaces.
Let
$H_1(2)\subset{\cal A}_3(2)$ be the inverse image of $H_1$ in ${\cal A}_3(2)$.
We denote the closure of $H_1$, resp. $H_1(2)$ in ${\cal A}_3^{\ast}$,
resp. ${\cal A}_3^{\ast}(2)$ by $\bar H_1$, resp. $\bar H_1(2)$.
\begin{proposition}\label{prop25}
The surface $\bar H_1(2)$ has $20$ components $\bar H^i_1(2)$ which are
equivalent under the group $\Gamma_3/\Gamma_3(2)\cong S_6$. Every component
$\bar H^i_1(2)$ is smooth and isomorphic to the product $X(2)\times X(2)$ of
modular curves of level $2$. Two different components ${\bar H}^i_1(2)$ and
${\bar H}^j_1(2)$ with $i\neq j$ do not intersect. The intersection of ${\bar
H}^i_1(2)$ with the boundary is transversal and consists of the curves
$\{\mbox{cusp}\}\times X(2)$ and $X(2)\times\{\mbox {cusp}\}$. The surface
${\bar
H}_1(2)$ is disjoint from the singularities of ${\cal A}^{\ast}_3(2)$.
\end{proposition}

\begin{Proof}
The number of components was computed in \cite[Theorem 3.2]{HSN}. Clearly the
intersection $H^i_1(2)$ of ${\bar H}^i_1(2)$ with ${\cal A}_3(2)$ is isomorphic
$X^{\circ}(2)\times X^{\circ}(2)$. To see that ${\bar H}^i_1(2)$ is still a
product one can proceed as in the proof of \cite[Proposition
(I.5.53)]{HKW1}. The
crucial point is to study the points in the boundary components $D(l)$, resp.
$E(h)$ which came from points with a non-trivial stabilizer. For $D(l)$ we did
this explicitly in the proof of Theorem (\ref{theo24}) which also gives us
immediately that ${\bar H}^i_1(2)$ and $D(l)$ intersect transversally. The
calculations for $E(h)$ are also straightforward although slightly more
cumbersome. Since this follows from \cite[Hilfssatz(II.5.23)]{Br} we will not
give the details here. At the same time we find that there are no points whose
stabilizer contains $\ZZ/2\times\ZZ/2$ and hence two different components
${\bar H}^i_1(2)$ and ${\bar H}^j_1(2)$  cannot meet. The singularities of
${\cal A}^{\ast}_3(2)$ are in the deepest points. None of these points have a
stabilizer which contains a group $\ZZ/2$ hence ${\bar H}_1(2)$ does not
contain
any of these points.
\hfill\end{Proof}

The curve $X(2)$ is rational and hence the above proposition shows that
${\bar H}_1(2)$ is a disjoint union of 20 quadrics. Since $X(2)$ has 3
cusps the
intersection of every such quadric with the boundary is a union of 6 lines
which
give a divisor of bidegree (3, 3). On the other hand the intersection of
${\bar H}_1$ with a boundary surface consists of 4 sections, namely the fixed
points of the Kummer involution. Altogether we have 15+15=30 boundary surfaces
and in this way we obtain a $(30_4, 20_6)$-configuration.\\

We now turn to the canonical bundle of ${\cal A}^{\ast}_3(2)$. Since the map
$\pi:\HH_2\rightarrow{\cal A}_3(2)$ is branched of order 2 along ${\bar
H}_1(2)$ it follows that
\begin{eqnarray}
K_{{\cal A}_3^{\ast}(2)} = 3L-D-\frac 12 H_1.
\end{eqnarray}
Here $L$ is the $\QQ$-line bundle of modular forms of weight 1 and $D$ denotes
the boundary.

\begin{lemma}\label{lem26}
For every component of ${\bar H}_1(2)$ the normal bundle has
bidegree $(-1,-1)$.
\end{lemma}
\begin{Proof}
Recall that for every component ${\bar H}^i_1(2)\cong X(2)\times X(2)$ is
a quadric and hence
$K_{{\bar H}^i_1(2)}$ has bidegree $(-2,-2)$. By adjunction
$$
K_{{\bar H}_1^i(2)} = \left(K_{{\cal A}_3^{\ast}(2)}+{\bar
H}^i_1(2)\right)|_{{\bar H}_1^i(2)}.
$$
Using formula (1) for $K_{{\cal A}_3^{\ast}(2)}$ we can deduce from this
that
\begin{eqnarray}
{\bar H}^i_1(2)|_{{\bar H}^i_1(2)}=2(K_{{\bar H}^i_1(2)}-3L+D)|_{{\bar
H}^i_1(2)}.
\end{eqnarray}
The line bundle $L$ restricts to $X(2)\times X(2)$ as $L_{X(2)}\boxtimes
L_{X(2)}$ which has bidegree
$(\frac 12, \frac 12)$. Since $X(2)$ has 3 cusps we had already seen that the
boundary intersects ${\bar H}^i_1(2)$ in a divisor of bidegree (3,3). Finally
$K_{{\bar H}^i_1(2)}$ has bidegree (-2, -2) and the claim follows immediately
by adding the bidegrees in the right hand side of formula (2).
\hfill\end{Proof}
\begin{theorem}\label{theo27}
There exists a commutative diagram
$$
\begin{CD}
{\tilde Z}@>f>> Z\\
@V{\tilde\pi}VV @VV\pi V\\
{\cal A}^{\ast}_3(2)@>f'>> Z'
\end{CD}
$$
where the vertical maps $\pi$ and ${\tilde\pi}$ are blow-ups of the
singularities and where $f$ and $f'$ contract each of the $20$ quadrics ${\bar
H}_1^i(2)$ to a line such that $Z$ is a smooth projective Calabi-Yau threefold.
\end{theorem}
\begin{Proof}
We first remark that the singularities are harmless: If $X$ is a singularity of
type $V_{\frac 13(1,1,1)}$ and $f:Y\rightarrow X$ is the blow-up of $X$ in the
singular point then $Y$ is smooth and $f^{\ast} K_X=K_Y$. Since the normal
bundle of every component ${\bar H}^i_1(2)$ has bidegree (1,1) we can contract
each set of rulings and obtain a manifold. The crucial point is to show that
this can be done in such a way that the resulting manifold is projective. We
shall do this using Cornalba's criterion \cite[Theorem 2]{Co}. Let $\Phi$ and
$\Psi$ be the rulings of a component ${\bar H}^i_1(2)$. We choose an ample
line bundle $A$ on ${\cal A}^{\ast}_3(2)$ and assume that $A.\Psi\ge
A.\Phi$. Since the boundary components cut out rulings on ${\bar H}^i_1(2)$ of
the form $\{\mbox{cusp}\}\times X(2)$ or $X(2)\times\{\mbox{cusp}\}$ we can
choose a component $D_k$ such that $D_k.\Phi=0$ and $D_k.\Psi=1$. Choose
$\alpha\gg 0$ such that $\alpha A+D_k$ is ample. Then we set
$$
s:=(\alpha A+D_k).\Phi=\alpha A.\Phi>0.
$$
With this choice of $s$ we find that
$$
\begin{array}{rcl}
(\alpha A+D_k+s{\bar H}^i_1(2)).\Phi & = & \alpha A.\Phi-s=0\\
(\alpha A+D_k+s{\bar H}^i_1(2)).\Psi & = & \alpha A.\Psi+1-s>0.
\end{array}
$$
It then follows from Cornalba's theorem that the manifold which is obtained by
contracting the rulings homologous to $\Phi$ is smooth projective. We can do
this for every component of ${\bar H}_1(2)$ separately.\\
The modular form $\Delta_1^3$ defines a canonical form
$$
\omega=\Delta_1^3 d \tau_1 \wedge d\tau_2\wedge d\tau_3
$$
which descends to ${\cal A}^{\ast}_3(2)$ and vanishes of order 1 along ${\bar
H}^i_1$. It has no other zeroes in ${\cal A}_3(2)$. Since the Fourier expansion
of $\Delta_1$ starts with $e^{2\pi i\tau_1/6}$, resp. $e^{2\pi i\tau_3/18}$
it follows that $\Delta^3_1$ vanishes of order 1 along the boundary $D$ and
hence $\omega$ has no zeroes along a boundary component. This proves that
$f'_{\ast}(\omega)$ defines a nowhere vanishing canonical form outside
$f'({\bar
H}_1(2))$. Since $f'({\bar H}_1(2))$ is a union of curves we have in fact
obtained a nowhere vanishing form on $Z'$, i.e. $K_{Z'}={\cal O}_{Z'}$ and
hence
also
$K_Z={\cal O}_Z$ since the only singularities present are of type $V_{\frac
13}(1, 1, 1)$. Finally
$q(Z)=h^1({\cal O}_Z)=0$ since this holds for every smooth projective model
of a
Siegel modular variety \cite{F2}.
\hfill\end{Proof}

\begin{remark}\label{rem28}
The Calabi-Yau manifold $Z$ is clearly birational to the Calabi-Yau model
constructed by Barth and Nieto in \cite{BN}, but it is not clear whether
the two
models are isomorphic. By recent results of Batyrev \cite{Bat} and
Huybrechts \cite{Huy} it follows, however, that the two models have the same
Betti numbers, in particular $e(Z)=80$.
\end{remark}

\section{The moduli space ${\cal A}_3(3)$}
We consider the group
$$
\on{Sp}^{(3)} (\Lambda_3, \ZZ)=\{g\in\on{Sp}(\Lambda_3, \ZZ); g\equiv {\bf 1}
\mbox{ mod } 3\}
$$
resp. the conjugated group
$$
\Gamma_3(3)=R_3^{-1}\on{Sp}(\Lambda_3, \ZZ) R_3
$$
where $R_3= \mbox{ diag } (1, 1, 1, 3)$. Then ${\cal A}_3(3)=\Gamma_3(3)
\backslash\HH_2$ is the moduli space of $(1,3)$-polarized abelian surfaces
with a level $3$ structure. (For the definition of generalized level $n$
structure see also \cite[chapter 8]{LB}.) Again we denote by ${\cal
A}_3^{\ast}(3)$ the Voronoi compactification of ${\cal A}_3(3)$.

\begin{theorem}\label{theo31}
{\rm(i)} The variety ${\cal A}_3^{\ast}(3)$ is Gorenstein and its singular
locus consists of
a finite number of points of type $V_{{\frac 13}(1, 1, 1)}$.\\
\noindent
{\rm(ii)} The variety ${\cal A}_3^{\ast}(3)$ is of general type.
\end{theorem}

\begin{Proof}
(i)  By \cite[Corollary 5.1.10]{LB} the group $\Gamma_3(3)$ acts freely on
$\HH_2$. To compute the singularities of ${\cal A}_3^{\ast}(3)$ we can then
proceed as in the proof of Theorem (\ref{theo24}). These arguments show
that for
$l=l_0$ and $l=l_{(0,1)}$ the group $P''_l(\Gamma_3(3))$ acts freely. The
analysis for the corank $2$ boundary components is also analogous. In this case
the lattice
$P'(h)$ is
$$
P'(h)=\left\{
\left(
\begin{array}{cc}
{\bf 1}_2  &  B\\
0          &  {\bf 1}_2
\end{array}
\right);\quad
B={^t B}, \quad B\in\left(
\begin{array}{cc}
3\ZZ  &  9\ZZ\\
9\ZZ  &  9\ZZ
\end{array}
\right)
\right\}
$$
and we find again a finite number of isolated singularities of type
$V_{\frac 13 (1,1, 1)}$.\\
(ii)  By (i) every modular form of weight $3k$
with respect to the group $\Gamma_3(3)$ which vanishes of order $k$ along
the boundary gives rise to a $k$-fold
differential form
on ${\cal A}_3^{\ast}(3)$. The form $\Delta^6_1$ is a modular form with respect
to $\Gamma_3(3)$ and vanishes of order $3$ along the boundary. Hence the space
$$
W_{9k}:=\Delta_1^{6k} M_{3k}(\Gamma_{1,3}(3)) \subset M_{9k}(\Gamma_{3}(3))
$$
is a space of forms of weight $9k$ vanishing of order $3k$ along the boundary.
Since the dimension of $W_{9k}$ grows as $\mbox{const}{\cdot} k^3$ it
follows that the
$3$-fold ${\cal A}_3^{\ast}(3)$ is of general type.
\hfill\end{Proof}

In this case we can again contract the components ${\bar H}^i_1(3)$ of the
Humbert surface of discriminant $1$ to curves to obtain a projective minimal
model of ${\cal A}_3^{\ast}(3)$.

\begin{theorem}\label{theo32}
The components ${\bar H}^i_1(3)$ of the Humbert surface ${\bar H}_1(3)$ of
discriminant $1$ can be contracted to curves such that the resulting variety
${\hat{\cal A}}_3(3)$ is projective and has finitely many isolated
singularities
of type $V_{\frac 13 (1, 1, 1)}$. The variety ${\hat{\cal A}}_3(3)$ is a
minimal
model of ${\cal A}_3^{\ast}(3)$.
\end{theorem}

\begin{Proof}
As before every component ${\bar H}_1^i(3)$ is isomorphic to the product
$X(3)\times X(3)$ of modular curves of level $3$, and hence to a quadric. Since
$\Gamma_3(3)$ acts freely on $\HH_2$ we find for the canonical bundle of ${\cal
A}_3^{\ast}(3)$ that
$$
K_{{\cal A}_3^{\ast}(3)}=3L-D.
$$
Using adjunction for the surfaces ${\bar H}^i_1(3)$ one obtains
$$
{\bar H}^i_1(3)|_{{\bar H}^i_1(3)}=K_{{\bar H}^i_1(3)}-(3L-D)|_{{\bar
H}^i_1(3)}.
$$
The line bundle $L_{X(3)}$ has degree $1$ and $X(3)$ has $4$ cusps. It follows
that the normal bundle of ${\bar H}_1^i(3)$ has bidegree $(-1, -1)$. We can now
argue as in the proof of Theorem \ref{theo27} and contract one set of rulings
for each of the quadrics ${\bar H}_1^i(3)$ to obtain a projective 3-fold
${\hat{\cal A}}_3(3)$. Since the singularities of ${\cal A}_3^{\ast}(3)$ are
disjoint from the Humbert surface ${\bar H}_1(3)$ the variety ${\hat{\cal
A}}_3(3)$ has the same singularities as ${\cal A}_3^{\ast}(3)$. In particular
${\hat{\cal A}}_3(3)$ is Gorenstein.

It remains to prove that $K_{{\hat{\cal A}}_3(3)}$ is nef. Let $f:{\cal
A}_3^{\ast}(3)\rightarrow{\hat{\cal A}}_3(3)$ be the contraction map. Then
$$
f^{\ast} K_{{\hat{\cal A}}_3(3)}=K_{{\cal A}_3^{\ast}(3)}-\sum\limits_i{\bar
H}_1^i(3)=K_{{\cal A}_3^{\ast}(3)}-{\bar H}_1(3)=:K'.
$$
We claim that $K'.C\ge 0$ for every curve $C$ in ${\cal A}_3^{\ast}(3)$. Since
$K'$ is trivial on ${\bar H}_1^i(3)$ we can assume that $C$ is not contained in
the Humbert surface of discriminant $1$. The form $\Delta^6_1$ is a modular
form (without a character) of weight $6$ with respect to the group
$\Gamma_3(3)$. It vanishes of order $6$ along the components ${\bar H}_1^i(3)$
of ${\bar H}_1(3)$ and of order $3$ along the boundary. Hence
$$
6L = 6{\bar H}_1(3) + 3D
$$
resp.
$$
D=2L-2{\bar H}_1(3).
$$
This implies that
$$
K_{{\hat{\cal A}}_3(3)}=L+2{\bar H}_1(3)
$$
resp.
$$
K'=L+{\bar H}_1(3).
$$
Since $L$ is nef (a multiple of $L$ is free) it follows that $K'.C\ge 0$ for
every curve $C$ not contained in ${\bar H}_1(3)$.
\end{Proof}

\section{Concluding remarks}
In this section we want to remark briefly on the geometric relevance of the
modular forms $\Delta_2$ and $\Delta_5$.\\

The modular form $\Delta_2$ is a modular form of weight $2$ with a character
of order $4$ with respect to the group $\Gamma_2$. The restriction of this
character to the congruence subgroup $\Gamma_2(4)$ is trivial. In particular
$\Delta_2$ is a weight $2$ cusp form with respect to $\Gamma_2(4)$ and argueing
along the same lines
as in the proof of Theorem \ref{theo31} this implies, under the assumption
that all singularities of ${\cal A}_2^{\ast}(4)$ are canonical, that
this space is of general type. The normal bundle of the
components ${\bar H}^i_1(4)$
of the Humbert surface of discriminant $1$ has bidegree
$(-2,-2)$ and the restriction of $K_{{\cal A}_2^{\ast}(4)}$ to these
components is trivial. Since $2L={\bar H}_1(4)-D$ it follows that
$K_{{\cal A}_2^{\ast}(4)}=L+{\bar H}_1(4)$ and, again modulo checking that all
singularities of ${\cal A}_2^{\ast}(4)$ are canonical, this shows that
the variety ${\cal A}_2^{\ast}(4)$ is a minimal model.\\

The modular form $\Delta_5$ has weight $5$ and a character of order $2$ with
respect to the integer symplectic group $\Gamma_1$. Since $\Delta_5$ vanishes
on the Humbert surface $H_1$ we find that $10L=2H_1$ on
${\cal A}_1$. The divisors $5L$
and $H_1$ differ by $2$-torsion (cf. \cite{GH2}).
On the moduli space ${\cal A}_1^{\ast}(n)$
we have the relation $10L=2{\bar H}_1(n)+nD$. Since
$K_{{\cal A}_1^{\ast}(n)}=3L-D$ it follows that
$$
K_{{\cal A}_1^{\ast}(n)}=(3-\frac{10}{n})L+\frac{2}{n}D.
$$
Note also that the normal bundle of the components ${\bar H}_1(n)$
is negative for $n\leq3$, trivial for $n=4$ and positive for $n\geq5$.
This can be used to show that $K_{{\cal A}_1^{\ast}(n)}$ is nef for $n=4$ and
ample for $n\geq5$ (see also \cite{H}, \cite{Bo}).

\vspace{1cm}

%
%
\bibliographystyle{alpha}

\noindent

Authors' addresses:

\bigskip

\parbox{7cm}
{Klaus Hulek\\
Institut f\"ur Mathematik\\
Universit\"at Hannover\\
D 30060 Hannover\\
Germany\\
E-Mail: hulek{\symbol{64}}math.uni-hannover.de}
\hfill
\parbox{7cm}
{Valery Gritsenko\\
St. Petersburg Branch of Steklov\\
Mathematical Institute\\
Fontanka 27, St.Petersburg \\
191011 Russia\\
e-mail:  gritsenk{\symbol{64}}kurims.kyoto-u.ac.jp}

\end{document}